\documentclass[12pt]{article}

\usepackage{amsmath}
\usepackage{amsfonts}
\usepackage{amssymb}

\newtheorem{Def}{Definition}[section]
\newtheorem{Prop}[Def]{Proposition}
\newtheorem{Thm}[Def]{Theorem}
\newtheorem{Lem}[Def]{Lemma}

\newtheorem{Cor}[Def]{Corollary}

\newcommand{\proof}{\noindent {\it Proof.} \ } 
\newcommand{\qed}{\hfill $\Box$}
\newcommand{\qedpar}{\qed \par \bigskip}
\newcommand{\R}{{\mathbb R}}
\newcommand{\C}{{\mathbb C}}
\newcommand{\Z}{{\mathbb Z}}

\newcommand{\LA}{{\mathfrak a}}
\newcommand{\LE}{{\mathfrak e}}
\newcommand{\LF}{{\mathfrak f}}
\newcommand{\LG}{{\mathfrak g}}

\newcommand{\LK}{{\mathfrak k}}

\newcommand{\LN}{{\mathfrak n}}
\newcommand{\LP}{{\mathfrak p}}
\newcommand{\LS}{{\mathfrak s}}

\newcommand{\sll}{\mathrm{sl}}

\newcommand{\ad}{\mbox{\rm ad}}

\newcommand{\aut}{\mbox{\rm Aut}}

\newcommand{\id}{\mbox{\rm id}}
\newcommand{\tr}{\mbox{\rm tr}}
 
\newcommand{\no}[1]{\mbox{\rm (#1)}} 
 
\newcommand{\inner}[2]{\langle #1 , #2 \rangle} 
\newcommand{\ric}{\mbox{\rm Ric}}

\numberwithin{equation}{section}

\title{\bf Noncompact homogeneous Einstein manifolds 
attached to graded Lie algebras} 
\author{Hiroshi Tamaru} 
\date{\small 
Department of Mathematics, Hiroshima University, \\ 
1-3-1 Kagamiyama, Higashi-Hiroshima 739-8526, JAPAN \\ 
(e-mail: tamaru@math.sci.hiroshima-u.ac.jp)
} 

\begin{document}

\maketitle 

\footnotetext{
Mathematics Subject Classification (2000): 
53C30, 
22E25, 
53C25, 
53C35, 
} 


\begin{abstract} 
In this paper, 
we study the nilradicals of parabolic subalgebras of semisimple Lie algebras 
and the natural one-dimensional solvable extensions of them. 
We investigate the structures, curvatures and Einstein conditions of 
the associated nilmanifolds and solvmanifolds. 
We show that our solvmanifold is Einstein if the nilradical is of two-step. 
New examples of Einstein solvmanifolds 
with three-step and four-step nilradicals are also given. 
\end{abstract} 

\section{Introduction} 

A Riemannian manifold is called a {\it solvmanifold} 
if it admits a transitive solvable group of isometries. 
Up to now, all known examples of noncompact homogeneous Einstein 
manifolds are solvmanifolds. 
The study of Einstein solvmanifolds is currently very active. 
In particular, Heber (\cite{H}) deeply studied Einstein solvmanifolds 
and obtained many essential results. 
One of his theorems states that 
the moduli space of Einstein solvmanifolds may have large dimension. 
Despite of that, 
we do not know many explicit examples of Einstein solvamnifolds; 
in particular, Einstein solvmanifolds whose nilradicals 
have the nilpotency greater than two. 
The purpose of this paper is to introduce and study a class of solvmanifolds, 
and provide new examples of Einstein spaces. 
Our class contains Einstein solvmanifolds whose nilradicals 
are of two, three and four-step. 
We need to know explicit examples 
for the further study of Einstein solvmanifolds. 

Here we briefly recall some known examples of 
noncompact homogeneous Einstein manifolds. 
One strategy for finding Einstein solvmanifolds is 
to consider solvmanifolds with good geometric structures. 
Typical examples are symmetric spaces of noncompact type. 
Other examples are noncompact homogeneous K\"ahler-Einstein manifolds, 
which are modeled by solvable normal $j$-algebras (see \cite{PS}), 
and quaternionic K\"ahler solvmanifolds (see \cite{A2}, \cite{C}). 
The other strategy is to consider 
solvable extensions of good nilmanifolds. 
Typical examples are Damek-Ricci spaces, 
which are the one-dimensional solvable extensions of generalized Heisenberg groups 
(\cite{Bo}, \cite{DR}, see also \cite{BTV}). 
Other examples are given by Mori (\cite{M}); 
he started from complex graded Lie algebras of second kind 
and considered the solvable extensions of the associated nilmanifolds. 
In fact, our construction is a generalization of his study. 
Other strategy for finding Einstein solvmanifolds 
can be found in recent articles, 
for example, see \cite{L-ann}, \cite{L} and \cite{GK}. 
The list of known examples of Einstein solvmanifolds can be found in \cite{L-dga}. 

In this paper we employ the second strategy. 
We take the nilradicals of parabolic subalgebras 
of semisimple Lie algebras. 
These nilpotent Lie algebras admit the natural inner products, 
and can be investigated in term of graded Lie algebras. 
In Section \ref{nil}, 
we mention some properties of graded Lie algebras 
and define our nilmanifolds. 
The structures and curvatures of our nilmanifolds will be studied 
in Sections \ref{sec-ricci} and \ref{sec-mean}. 
In Section \ref{sec-solv}, 
we study curvature property of one-dimensional solvable extensions 
of our nilmanifolds. 
We find that there is a natural solvable extension, 
that is, only one candidate to be Einstein. 
Curvatures of the natural solvable extension is studied in 
Section \ref{sec-naturalextension}. 
We see in Section \ref{sec-two} that 
our solvmanifold is Einstein if the nilradical is of two-step. 
This provides many new examples of Einstein solvmanifolds 
and some of them have "unusual" eigenvalue types. 
Other new examples of Einstein solvmanifolds 
with three-step and four-step nilradicals 
will be given in Section \ref{higher-step}. 
We observe in Section \ref{sec-lowrank} that, 
if $\LG$ is small, then every attached solvmanifolds are Einstein.

Note that, if a parabolic subalgebra is minimal, 
then its nilradical coincides with the nilpotent part 
of the Iwasawa decomposition. 
In this case the constructed solvmanifold 
coincides with the rank one reduction of a symmetric space of noncompact type. 
Therefore our class of solvmanifolds essentially 
contains noncompact symmetric spaces. 
The structures of the rank one reductions of noncompact symmetric spaces 
will be mentioned in Section \ref{sec-symm}. 

The author would like to thank 
Professor Soji Kaneyuki, Professor Hiroshi Asano and the referee. 
Theorem \ref{mcv-3} was pointed out by Professor Asano and the referee. 
Theorem \ref{ric-s-1} was suggested by the referee. 
Their advices are valuable to refine the first version of this manuscript. 
This work was partially supported by 
Grant-in-Aid for Young Scientists (B) 14740049 and 17740039, 
The Ministry of Education, Culture, Sports, Science and Technology, Japan. 

\section{Preliminaries} 
\label{pre}

In this section we recall some notations and basic properties of 
Einstein solvmanifolds. 
We refer to \cite{W}, \cite{H}. 

Let $\LS$ be a solvable Lie algebra and 
$\inner{}{}$ be a positive definite inner product on $\LS$. 
Throughout this paper, by a metric Lie algebra $(\LS, \inner{}{})$ 
we also mean the corresponding solvmanifold, 
that is, the simply-connected Lie group $S$ 
equipped with the induced left-invariant Riemannian metric. 
Put $\LN := [\LS,\LS]$ and $\LA := \LN^{\perp}$. 

\begin{Def} 
{\rm 
A solvmanifold $(\LS, \inner{}{})$ is called {\it standard} 
if $\LA$ is abelian. 
A standard solvmanifold $(\LS, \inner{}{})$ is said to be of 
{\it Iwasawa-type} if \ 
(i) $\ad_A$ is symmetric for every $A \in \LA$, \ and 
(ii) there exists $A_0 \in \LA$ such that 
every eigenvalue of $\ad_{A_0}$ on $\LN$ is a positive real number. 
} 
\end{Def} 

To calculate the Ricci curvatures of our solvmanifolds, 
we need the following particular vectors. 

\begin{Def} 
{\rm 
For a solvmanifold $(\LS, \inner{}{})$, 
we call $H_0 \in \LA$ the {\it mean curvature vector} if 
$\inner{H_0}{A} = \tr (\ad_A)$ for every $A \in \LA$. 
} 
\end{Def} 

Note that $H_0$ is the mean curvature vector 
for the embedding of $N$ into $S$, 
where $N$ denotes the Lie subgroup of $S$ with Lie algebra $\LN$. 

\begin{Prop}[\cite{W}] 
\label{wolter} 
Let $(\LS, \inner{}{})$ be a solvmanifold of Iwasawa type 
and $H_0$ be the mean curvature vector. 
Denote by $\ric$ and $\ric^{\LN}$ the Ricci curvatures of 
$\LS$ and $\LN$, respectively. 
Then, 
\begin{enumerate} 
\item[\no{1}] \ 
$\ric(A,A') = - \tr (\ad_A) \circ (\ad_{A'})$ 
for $A,A' \in \LA$, 
\item[\no{2}] \ 
$\ric(A,X) = 0$, 
\item[\no{3}] \ 
$\ric(X) = \ric^{\LN}(X) - (\tr(\ad_{H_0}) / |H_0|^2) 
\cdot \ad_{H_0}X$ \ for $X \in \LN$.
\end{enumerate} 
\end{Prop} 

Proposition \ref{wolter} leads that, 
if $(\LS, \inner{}{})$ is Einstein 
then so is $(\R H_0 + \LN, \inner{}{})$. 
For a solvmanifold $(\LS, \inner{}{})$, 
we call $\dim \LA$ the {\it algebraic rank} 
and $(\R H_0 + \LN, \inner{}{})$ the {\it rank one reduction}. 
Note that $H_0$ is also the mean curvature vector of the solvmanifold 
$(\R H_0 + \LN, \inner{}{})$. 

\begin{Prop}[\cite{A}] 
\label{alek} 
The Ricci curvature of a nilmanifold $(\LN, \inner{}{})$ is given by 
\[ \ric^{\LN} = 
(1/4) \sum \ad_{E_i} \circ (\ad_{E_i})^{\ast} 
- (1/2) \sum (\ad_{E_i})^{\ast} \circ \ad_{E_i} , \] 
where $\{ E_i \}$ is an orthonormal basis of $\LN$ 
and $(\ad_{E_i})^{\ast}$ denotes the dual of the adjoint operator 
defined by 
$\inner{(\ad_{E_i})^{\ast} X}{Y} = \inner{X}{\ad_{E_i} Y}$. 
\end{Prop} 

Propositions \ref{wolter} and \ref{alek} state that 
the Ricci curvatures of solvmanifolds of Iwasawa-type 
can be calculated in terms of the mean curvature vectors, 
adjoint operators and their dual operators. 

\begin{Prop}[\cite{H}] 
For a standard Einstein solvmanifold, 
there exists $\lambda >0$ such that 
the real part of every eigenvalue of $\ad_{\lambda H_0}$ is a positive integer. 
\end{Prop} 

Let $\mu_1 < \mu_2 < \cdots < \mu_m$ be the real parts of eigenvalues of 
$\ad_{\lambda H_0}$ without common divisor, 
and $d_1, d_2, \ldots, d_m$ the corresponding multiplicities. 
We call 
$(\mu_1 < \cdots < \mu_m ; d_1, \ldots, d_m)$ 
the {\it eigenvalue type} of an Einstein solvmanifold. 
For example, the eigenvalue types of Damek-Ricci spaces 
are $(1<2;d_1,d_2)$. 

\section{Nilmanifolds attached to graded Lie algebras} 
\label{nil}

In this section we briefly recall some properties of graded Lie algebras 
and define the nilmanifolds associated to them. 
We refer to \cite{K}, \cite{T4} for graded Lie algebras. 

Let $\LG$ be a semisimple Lie algebra. 
A decomposition of $\LG$ into subspaces, 
\[ \LG = \sum_{k \in \Z} \LG_k \] 
is called a {\it gradation} if $[\LG_k, \LG_l] \subset \LG_{k+l}$ 
holds for every $k$ and $l$. 
A Lie algebra with gradation is called a {\it graded Lie algebra}. 
By definition, $\LN := \sum_{k>0} \LG_k$ is a nilpotent subalgebra. 

In the theory of graded Lie algebras, 
it is fundamental that there exist the characteristic element 
$Z \in \LG$ and a grade-reversing Cartan involution $\sigma$. 
An element $Z$ is called the {\it characteristic element} if 
$\LG_k = \{ X \in \LG \mid [Z,X] = k \cdot X \}$ for every $k$. 
A Cartan involution $\sigma$ is said to be {\it grade-reversing} 
if $\sigma(\LG_k) = \LG_{-k}$ for every $k$. 
Note that the characteristic element is unique. 
It is easy to show that $\sigma$ is grade-reversing 
if and only if $\sigma(Z)=-Z$. 
A Cartan involution $\sigma$ defines the natural inner product 
$B_{\sigma}$ on $\LG$ by $B_{\sigma}(X,Y) := -B(X,\sigma(Y))$, 
where $B$ denotes the Killing form of $\LG$. 
Now we can define our nilmanifolds: 

\begin{Def} 
{\rm 
Let $\LG = \sum \LG_k$ be a graded Lie algebra. 
The nilpotent Lie algebra $\LN := \sum_{k>0} \LG_k$ 
endowed with the inner product $B_{\sigma}$ is called 
the {\it nilmanifold attached to the graded Lie algebra}, 
or shortly, {\it nilmanifold attached to GLA}. 
} 
\end{Def} 

Note that $\sum_{k \geq 0} \LG_k$ is a parabolic subalgebra of $\LG$ 
and $\LN$ is the nilradical of it. 
Conversely, it is known that every parabolic subalgebra of $\LG$ 
can be obtained in this way. 
Therefore, we are considering the nilradicals of 
all parabolic subalgebras of all semisimple Lie algebras. 

We also note that our nilpotent Lie algebras $\LN$ are 
related to $R$-spaces. 
Let $\LG = \sum \LG_k$ be a graded Lie algebra, 
$G$ the centerless Lie group with Lie algebra $\LG$, 
and $P$ the subgroup of $G$ with Lie algebra $\sum_{k \leq 0} \LG_k$. 
The quotient space $G/P$ is an $R$-space 
(it is sometimes called a real flag manifold). 
By definition, $\LN$ can be naturally identified 
with the tangent spaces of $G/P$. 
Therefore, there is a correspondence between 
the class of our nilmanifolds and the class of $R$-spaces 
(but this is not one-to-one). 

A graded Lie algebra $\LG = \sum \LG_k$ 
is said to be of {\it $\nu$-th kind} if 
$\LG_{\nu} \neq 0$ and $\LG_k =0$ for every $k > \nu$. 
The existence of a grade-reversing Cartan involution 
implies that $\dim \LG_k = \dim \LG_{-k}$, 
therefore, a $\nu$-th kind gradation satisfies 
$\LG_k =0$ for every $|k| > \nu$. 
Note that, for a graded Lie algebra of $\nu$-th kind, 
the nilpotent Lie algebra $\LN = \sum_{k>0} \LG_k$ 
is not of $\nu$-step in general, but at most $\nu$-step nilpotent. 
We will explain this later. 


To state the classification of gradations, 
we now start from a semisimple symmetric pair of noncompact type 
$(\LG, \LK)$ with Cartan involution $\sigma$. 
Take the Cartan decomposition $\LG = \LK + \LP$ 
with respect to $\sigma$ and a maximal abelian subspace $\LA$ in $\LP$. 
Let $\Delta$ be the root system of $\LG$ with respect to $\LA$. 
Recall that we say $\alpha \in \LA^{\ast} - \{ 0 \}$ is a {\it root} if 
\[ 0 \neq \LG_{\alpha} := \{ X \in \LG \mid [H,X] = \alpha(H) \cdot X \
\mbox{for every $H \in \LA$} \} . \]
Note that $\LG$ can be decomposed into the direct sum of root spaces: 
$\LG = \LG_0 + \sum_{\alpha \in \Delta} \LG_{\alpha}$. 
Let $\Pi := \{ \alpha_1, \ldots , \alpha_r \}$ be 
a set of simple roots of $\Delta$ and define 
\begin{eqnarray*} 
C & := & \{ H \in \LA \mid \alpha_i(H) \in \Z_{\geq 0} \ 
\mbox{for every $\alpha_i \in \Pi$} \} \\ 
& = & \{ c_1 H^1 + \cdots + c_r H^r \mid 
c_1, \ldots, c_r \in \Z_{\geq 0} \} , 
\end{eqnarray*} 
where $\{ H^1, \ldots, H^r \}$ is the dual basis of $\Pi$, 
that is, $\alpha_i(H^j) = \delta_{ij}$. 
Then, every element $Z \in C$ gives a gradation, by putting 
\[ \LG_k := \sum_{\alpha(Z)=k} \LG_{\alpha} . \] 
On the other hand, every gradation of $\LG$ can be obtained 
in this way up to conjugation (see \cite{KA}, \cite{T4}). 
We summarize this as follows. 

\begin{Thm}[\cite{KA}, \cite{T4}] 
There is a bijective correspondence between 
the set of the isomorphism classes of gradations on $\LG$ 
and the set of $\aut(DD)$-orbits of $C$, 
where $\aut(DD)$ denotes the automorphism group 
of the Dynkin diagram of $\Delta$. 
\end{Thm} 

It is easy to see that the gradation given by $Z \in C$ 
is of $\nu$-th kind if and only if $\nu = \widetilde{\alpha}(Z)$, 
where $\widetilde{\alpha}$ denotes the highest root. 
Obviously $\LN = \sum_{k>0} \LG_k$ is at most $\nu$-step nilpotent. 
We will mention here the condition for $\LN$ 
to be $\nu$-step nilpotent. 

A graded Lie algebra $\LG = \sum \LG_k$ is said to be 
{\it of type $\alpha_0$} 
if $\LG_1$ generates $\sum_{k>0} \LG_k$. 
This notation was introduced by Kaneyuki and Asano (\cite{KA}), 
and they showed that a graded Lie algebra $\LG = \sum \LG_k$ is 
of type $\alpha_0$ if and only if 
$Z \in \{ c_1 H^1 + \cdots + c_r H^r \mid 
c_1, \ldots, c_r \in \{ 0,1 \} \}$ up to conjugation. 
This property is related to the structure 
of the attached nilmanifold. 

\begin{Prop} 
Let $(\LN, B_{\sigma})$ be the nilmanifold attached to 
a graded Lie algebra $\LG = \sum \LG_k$ of $\nu$-th kind. 
Then, 
\begin{enumerate} 
\item[\no{1}] \ 
there exists a type $\alpha_0$ gradation $\LG = \sum \LG'_k$ 
which satisfies $\LN = \sum_{k>0} \LG'_k$, 
\item[\no{2}] \ 
$\LN$ is of $\nu$-step nilpotent 
if and only if the gradation is of type $\alpha_0$. 
\end{enumerate} 
\end{Prop} 

\proof 
(1): Denote the characteristic element of the gradation 
$\LG = \sum \LG_k$ by 
$Z = c_{i_1} H^{i_1} + \cdots + c_{i_s} H^{i_s}$ 
with $c_{i_1}, \ldots, c_{i_s} >0$. 
Then the gradation given by 
$Z' = H^{i_1} + \cdots + H^{i_s}$ 
satisfies the condition. 
(2): Assume $\LG = \sum \LG_k$ is of type $\alpha_0$. 
Since $\LG_1$ generates $\LG_{\nu}$, 
the lower central series satisfies $C^{\nu}(\LN) \neq 0$. 
Thus $\LN$ is of $\nu$-step nilpotent. 
To show the converse, assume $\LG = \sum \LG_k$ is not 
of type $\alpha_0$. 
In this case one has $Z = c_{i_1}H^{i_1} + \cdots + c_{i_s}H^{i_s}$ 
and some coefficients are greater than $1$. 
As in the proof for (1), 
the gradation $\LG = \sum \LG'_k$ given by $H^{i_1} + \cdots + H^{i_s}$ 
is of type $\alpha_0$ and satisfies $\LN = \sum_{k>0} \LG'_k$. 
The assumption implies that 
$\LG = \sum \LG'_k$ is of $\nu'$-kind with $\nu' < \nu$. 
Therefore $\LN$ is $\nu'$-step nilpotent and not $\nu$-step nilpotent. 
\qedpar 

This proposition means that, 
to study the nilmanifold attached to GLA, 
it is enough to consider gradations of type $\alpha_0$. 

\section{Ricci curvatures of the attached nilmanifolds} 
\label{sec-ricci}

In this section we calculate the Ricci curvatures 
of the nilmanifolds attached to GLA. 
The Ricci curvatures can be written by using 
the bracket products and the grade-reversing Cartan involutions. 

\begin{Lem} 
\label{lem-n} 
The nilmanifolds $(\LN, B_{\sigma})$ attached to GLA satisfy that 
\begin{enumerate} 
\item[\no{1}] \ 
$B_{\sigma}([U,V],W) = -B_{\sigma}(V,[\sigma U,W])$ 
\ for every $U,V,W \in \LN$, 
\item[\no{2}] \ 
$(\ad_U)^{\ast} V = [V,\sigma U]_{\LN}$ 
\ for every $U,V \in \LN$. 
\end{enumerate} 
\end{Lem} 

These formulae can be checked easily 
but quite useful for calculating the Ricci curvatures. 

\begin{Def} 
{\rm 
Let $\{ E^{(k)}_i \}$ be an orthonormal basis of $\LG_k$. 
We call 
$Z_k := \sum [\sigma(E^{(k)}_i),E^{(k)}_i]$ 
the {\it $k$-th mean curvature vector}.} 
\end{Def} 

It is easy to show that $Z_k$ is independent of the choice of 
orthonormal basis. 
Furthermore, Lemma \ref{lem-n} leads that 
the vector $H_0 := Z_1 + \cdots + Z_{\nu}$ 
coincides with the mean curvature vector of the 
solvmanifold $(\LA + \LN, B_{\sigma})$. 



\begin{Lem} 
\label{ricci-n} 
Let $\{ E^{(k)}_i \}$ be an orthonormal basis of $\LG_k$ and 
define the operators $A_k$ and $B_k$ on $\LN$ by 
\[ 
A_k(U) := \sum [E^{(k)}_i,[\sigma E^{(k)}_i,U]] , \quad
B_k(U) := \sum [\sigma E^{(k)}_i,[E^{(k)}_i,U]] . 
\] 
Then, for every $U_l \in \LG_l$, we have 
\begin{enumerate} 
\item[\no{1}] \ 
$A_k(U_l) = B_k(U_l) - [Z_k,U_l]$, 
\item[\no{2}] \ 
$B_k(U_l) = A_{k+l}(U_l) = B_{k+l}(U_l) - [Z_{k+l},U_l]$. 
\end{enumerate} 
\end{Lem} 

\proof 
(1) is a direct consequence of the Jacobi identity. 
To show (2), we need an orthonormal basis 
$\{ E^{(k+l)}_j \}$ of $\LG_{k+l}$. 
It follows from $[E^{(k)}_i,U_l] \in \LG_{k+l}$ that 
\begin{eqnarray*} 
B_m(U_l) 
& = & \sum_{i} [\sigma E^{(k)}_i, [E^{(k)}_i,U_l]] \\
& = & \sum_{i} \sum_{j} 
[\sigma E_i, \inner{[E^{(k)}_i,U_l]}{E^{(k+l)}_j} E^{(k+l)}_j] \\ 
& = & 
\sum_{i} \sum_{j} 
[\sigma(\inner{E^{(k)}_i}{[\sigma U_l,E^{(k+l)}_j]} E^{(k)}_i), 
E^{(k+l)}_j] \\
& = & \sum_{j} 
[\sigma([\sigma U_l,E^{(k+l)}_j]), E^{(k+l)}_j] \\
& = & A_{l+m} (U_l) . 
\end{eqnarray*} 
This concludes the lemma. 
\qedpar 

By definition, $B_k(U_l)=0$ if $k+l > \nu$. 
Therefore one can determine $A_k$ and $B_k$ 
by these recursion formulae. 
These operators play important roles for calculating the Ricci curvatures 
of our nilmanifolds. 

\begin{Thm} 
\label{ric-n-1} 
Let $(\LN, B_{\sigma})$ be the nilmanifold attached to GLA 
and $\ric^{\LN}$ be the Ricci curvature. 
Then, for every $U_l \in \LG_l$, one has 
\begin{enumerate} 
\item[\no{1}] \ 
$\displaystyle \ric^{\LN}(U_l) = -(1/4) \sum_{m<l} A_m(U_l)
+(1/2) \sum_{m>l} A_m(U_l)$, 
\item[\no{2}] \ 
$\ric^{\LN}(U_1) = (1/2) \sum (1-k) [Z_k,U_1]$, 
\item[\no{3}] \ 
$\ric^{\LN}(U_2) = (1/4) \sum (2-k) [Z_k,U_2]$. 
\end{enumerate} 
\end{Thm} 

\proof 
Let $\{ E_i \}$ be an orthonormal basis of $\LG_k$. 
Direct calculations lead that 
\[ 
\sum \ad_{E_i} \circ (\ad_{E_i})^{\ast} (U_l)
= - \sum [E_i, [\sigma E_i,U_l]_{\LN}] =
\left\{ 
\begin{array}{ll} 
-A_m & \mbox{for $m<l$} , \\
0 & \mbox{for $m \geq l$} .
\end{array} 
\right. 
\] 
Furthermore, one can see that 
\[ 
\sum (\ad_{E_i})^{\ast} \circ \ad_{E_i} (U_l)
= -B_m(U_l) = -A_{l+m}(U_l) . 
\] 
Thus Proposition \ref{alek} concludes (1). 
One can calculate $A_k(U_1)$ and $A_k(U_2)$ 
by using Lemma \ref{ricci-n} as follows: 
\begin{eqnarray*} 
A_k(U_1) & = & B_k(U_1) - [Z_k,U_1] \\
& = & (1/2) B_{k+1}(U_1) - [Z_k,U_1] - [Z_{k+1},U_1] \\
& = & \cdots \\ 
& = & - [Z_k + Z_{k+1} + \cdots + Z_{\nu},U_1] , \\ 
A_k(U_2) 
& = & -[Z_k, U_2] -[Z_{k+2}, U_2] -[Z_{k+4}, U_2] - \cdots . 
\end{eqnarray*} 
By substituting these formulae for the equation (1), 
we conclude (2) and (3). 
\qed 

\section{Mean curvature vectors} 
\label{sec-mean}

As we saw in the previous section, 
one needs to know the mean curvature vectors $Z_1, \ldots, Z_{\nu}$ 
for calculating the Ricci curvature of the nilmanifolds attached to GLA. 
In this section we investigate these vectors. 

We call $H_{\alpha} \in \LA$ the {\it root vector} of a root $\alpha$ 
if $B_{\sigma}(H_{\alpha},H) = \alpha(H)$ for every $H \in \LA$. 
Our root vectors are the negative of the usual ones, 
since $B_{\sigma} = -B$ on $\LA \times \LA$. 
We employ this definition for convenience. 
Let $\Delta_k := \{ \alpha \in \Delta \mid \alpha(Z)=k \}$. 

\begin{Lem} 
\label{4-1} 
The $k$-th mean curvature vector $Z_k$ satisfies 
\no{1} $Z_k = \sum_{\alpha \in \Delta_k} (\dim \LG_{\alpha}) H_{\alpha}$, 
and \no{2} $B_{\sigma}(Z_k,Z) = k \dim \LG_k$. 
\end{Lem} 

\proof 
Let $E \in \LG_{\alpha}$ be a unit vector. 
Then, Lemma \ref{lem-n} follows that $H_{\alpha} = [\sigma E,E]$, 
which concludes (1). 
To show (2) we need $\{ E_i \}$, an orthonormal basis of $\LG_k$. 
One can see that 
\begin{eqnarray*} 
B_{\sigma}(Z_k,Z)
& = & B_{\sigma}(\sum [\sigma E_i,E_i], Z)
\ = \ \sum B_{\sigma}(E_i,[Z,E_i]) \\ 
& = & \sum B_{\sigma}(E_i, k E_i) 
\ = \ k \dim \LG_k .  
\end{eqnarray*} 
This concludes (2). 
\qedpar 


We now determine the mean curvature vectors of our nilmanifolds. 
Recall that we can assume, without loss of generality, 
that $Z$ is a linear combination of $\{ H^i \}$, 
the dual basis of the simple roots $\Pi = \{ \alpha_i \}$. 
We also recall that the Weyl group is generated by $\{ s_1, \ldots, s_r \}$, 
where $s_j$ denotes the reflection of $\LA$ with respect to 
the hyperplane $\{ H \in \LA \mid \alpha_j(H)=0 \}$. 


\begin{Thm} 
\label{mcv-2} 
Let $(\LN, B_{\sigma})$ be the nilmanifold attached to GLA 
and assume that the characteristic element is 
$Z = c_{i_1} H^{i_1} + c_{i_2} H^{i_2} + \cdots + c_{i_s} H^{i_s}$ 
with $c_{i_1}, \ldots, c_{i_s} >0$. 
Then, 
\begin{enumerate}
\item[\no{1}] \ 
$Z_k = z_{i_1} H^{i_1} + \cdots + z_{i_s} H^{i_s}$ with 
$z_{i_j} = \sum_{\beta \in \Delta_k} 
(\dim \LG_{\beta}) B_{\sigma}(H_{\alpha_{i_j}},H_{\beta})$, 
for every $k$, and 
\item[\no{2}] \ 
$H_0 = h_{i_1} H^{i_1} + \cdots + h_{i_s} H^{i_s}$ with 
$h_{i_j} = \sum_{\beta \in \Delta^{i_j}} 
(\dim \LG_{\beta}) B_{\sigma}(H_{\alpha_{i_j}},H_{\beta})$, 
where $\Delta^{i_j} := \{ \beta \in \Delta \mid \beta(Z) = \beta(H^{i_j}) > 0 \}$. 
Furthermore, $h_{i_j} >0$. 
\end{enumerate}
\end{Thm} 

\proof 
(1): \ 
First of all we show that 
$Z_k \in V := \R H^{i_1} + \R H^{i_2} + \cdots + \R H^{i_s}$. 
Since the orthogonal complement of $V$ in $\LA$ 
is spanned by root vectors 
$\{ H_{\alpha_j} \mid j \neq i_1, \ldots, i_s \}$, 
it is sufficient to show that 
$B_{\sigma}(Z_k, H_{\alpha_j}) = 0$ for every $j \neq i_1, \ldots, i_s$. 
By the definition of Weyl group, 
this condition is equivalent to $s_j(Z_k)=Z_k$. 
Take $j \neq i_1, \ldots, i_s$. 
Because $s_j$ acts trivially on $V$, 
one has $s_j(Z)=Z$ and hence $s_j(\Delta_k)=\Delta_k$. 
This and Lemma \ref{4-1} imply that 
\[ 
s_j(Z_k) 
= s_j(\sum_{\beta \in \Delta_k} (\dim \LG_{\beta}) H_{\beta})
= \sum_{\beta \in \Delta_k} (\dim \LG_{\beta}) s_j(H_{\beta})
= \sum_{\beta \in \Delta_k} (\dim \LG_{s_j(\beta)}) H_{s_j(\beta)}
= Z_k . 
\] 
Note that the multiplicities of roots, $\dim \LG_{\beta}$, 
are invariant under the action of the Weyl group. 
Therefore we have $B_{\sigma}(Z_k, H_{\alpha_j}) = 0$, 
and hence $Z_k \in V$. 
Each coefficient of $Z_k = z_{i_1} H^{i_1} + \cdots + z_{i_s} H^{i_s}$ 
can be determined by 
\[ z_{i_j} = \alpha_{i_j}(Z_k) 
= \alpha_{i_j}(\sum_{\beta \in \Delta_k} 
(\dim \LG_{\beta}) H_{\beta}) 
= \sum_{\beta \in \Delta_k} 
(\dim \LG_{\beta}) B_{\sigma}(H_{\alpha_{i_j}},H_{\beta}) . \] 

(2): \ 
Let $\Delta' := \{ \beta \in \Delta \mid \beta(Z)>0 \}$. 
Because of (1), 
one can write $H_0 = h_{i_1} H^{i_1} + \cdots + h_{i_s} H^{i_s}$ 
with 
$h_{i_j} = \sum_{\beta \in \Delta'} (\dim \LG_{\beta}) 
B_{\sigma}(H_{\alpha_{i_j}},H_{\beta})$. 
We here claim that $s_{i_j}$ preserves $\Delta' - \Delta^{i_j}$. 
For every $\beta \in \Delta' - \Delta^{i_j}$, 
one can write $\beta = \sum n_k \alpha_k$ 
so that $n_{i_t} > 0$ for some $t \neq j$. 
The formula of the reflection 
\[ s_{i_j} (\beta) = \beta - 2 
\frac{\inner{\alpha_{i_j}}{\beta}}{\inner{\beta}{\beta}} \alpha_{i_j} \] 
leads that the $\alpha_{i_t}$-coefficient of $s_{i_j}(\beta)$ is also $n_{i_t}$. 
Hence, $s_{i_j}(\beta)(Z) \geq n_{i_t} >0$ and 
$s_{i_j}(\beta)(Z) \geq n_{i_t} + \beta(H^{i_j}) > \beta(H^{i_j})$. 
This completes the claim. 
Thus, the similar argument as (1) leads that 
\[ \sum_{\beta \in \Delta' - \Delta^{i_j}} (\dim \LG_{\beta}) 
B_{\sigma}(H_{\alpha_{i_j}},H_{\beta}) =0, \] 
which concludes the first part of (2). 
Now we show that $h_{i_j}>0$. 
There might exist $\beta \in \Delta^{i_j}$ such that 
$B_{\sigma}(H_{\alpha_{i_j}},H_{\beta}) <0$. 
But, for such $\beta$, one has $s_{i_j}(\beta) \in \Delta^{i_j}$ and 
$B_{\sigma}(H_{\alpha_{i_j}},s_{i_j}(H_{\beta})) = 
-B_{\sigma}(H_{\alpha_{i_j}},H_{\beta})$. 
Therefore every negative term cancels with another term. 
Hence, we got 
$h_{i_j} = \sum (\dim \LG_{\beta}) B_{\sigma}(H_{\alpha_{i_j}},H_{\beta})$, 
where $\beta$ run through 
$\{ \beta \in \Delta^{i_j} \mid 
B_{\sigma}(H_{\alpha_{i_j}},H_{\beta}) > 0, \ 
s_{i_j}(\beta) \not\in \Delta^{i_j} \}$. 
Since this set contains $\alpha_{i_j}$, one has $h_{i_j}>0$. 
\qedpar 

Next we show a remarkable property of mean curvature vectors, 
which is quite useful for calculating the Ricci curvatures 
of our nilmanifolds and solvmanifolds. 

\begin{Thm}
\label{mcv-3} 
Let $(\LN, B_{\sigma})$ be the nilmanifold attached to GLA 
with the characteristic element $Z$. 
Then the mean curvature vectors $Z_k$ satisfy 
\[ 
2 \sum_{k>0} kZ_k = Z . 
\]
\end{Thm}

\proof 
It is enough to show that $\inner{2 \sum kZ_k}{H} = \inner{Z}{H}$ 
for every $H \in \LA$. 
Let $E_i$ be an orthonormal basis of $\LG$. 
By the definition of the Killing form, one has 
\begin{eqnarray*}
\inner{Z}{H} 
= \tr (\ad_{Z} \circ \ad_{H}) 
= \sum_i \inner{[Z,[H,E_i]]}{E_i} . 
\end{eqnarray*}
One can assume that each $E_i$ is contained in the root spaces. 
Denote by $E_{\alpha} \in \LG_{\alpha}$ an unit vector. 
Thus we have 
\begin{eqnarray*}
\inner{Z}{H} 
& = & 
\sum_{\alpha \in \Delta} 
\inner{[Z,[H,E_{\alpha}]]}{E_{\alpha}} \dim \LG_{\alpha} \\ 
& = & 
\sum_{\alpha \in \Delta} \alpha(Z) \alpha(H) \dim \LG_{\alpha} . 
\end{eqnarray*}
On the other hand, Lemma \ref{4-1} (1) leads that 
\begin{eqnarray*}
2 \sum kZ_k
= 
2 \sum_{\alpha \in \Delta^{+}} \alpha(Z) (\dim \LG_{\alpha}) H_{\alpha} 
= 
\sum_{\alpha \in \Delta} \alpha(Z) (\dim \LG_{\alpha}) H_{\alpha} . 
\end{eqnarray*}
Therefore one can easily see that 
\[ 
\inner{2 \sum kZ_k}{H} = 
\sum_{\alpha} \alpha(Z) (\dim \LG_{\alpha}) \alpha(H), 
\]
which completes the proof. 
\qed 

\bigskip 

Combining this theorem with Theorem \ref{ric-n-1}, 
the Ricci curvatures of our nilmanifolds can be written easier. 

\begin{Cor}
\label{mcv-4}
For the nilmanifold $(\LN, B_{\sigma})$ attached to GLA, we have 
\begin{enumerate} 
\item[\no{1}] \ 
$\ric^{\LN}(U_1) = -(1/4) U_1 + (1/2) \sum [Z_k,U_1]$ 
\quad for $U_1 \in \LG_1$, 
\item[\no{2}] \ 
$\ric^{\LN}(U_2) = -(1/4) U_2 + (1/2) \sum [Z_k,U_2]$  
\quad for $U_2 \in \LG_2$. 
\end{enumerate}
\end{Cor}

%
%
%

\section{Curvature properties of the solvable extensions} 
\label{sec-solv}

We now consider one-dimensional solvable extensions of 
the nilmanifolds $(\LN, B_{\sigma})$ attached to GLA, 
and investigate their curvature properties and Einstein conditions. 
We consider the following solvable extensions: 
\[ 
\LS := \R H + \LN, \qquad 
\inner{}{}^c := c B_{\sigma} |_{\R H \times \R H}
+ B_{\sigma} |_{\LN \times \LN} , 
\] 
where $c>0$, $H \in \LA$ and $\ad_H |_{\LN}$ has only positive eigenvalues. 

A solvmanifold $(\LS, \inner{}{})$ is called a 
{\it Carnot space} if $\dim [\LS,\LS]^{\perp} = 1$ and 
there exists $A_0 \in [\LS,\LS]^{\perp}$ such that 
every eigenvalue of $\ad_{A_0}$ on $[\LS,\LS]$ is a positive integer (see \cite{N}). 
We call a Carnot space is of {\it $(k+1)$-step} if 
$k$ is the maximum of the eigenvalues of $\ad_{A_0}$. 

\begin{Prop}
The above solvmanifold $(\LS = \R H + \LN, \inner{}{}^c)$ is of Iwasawa-type. 
This solvmanifold is a Carnot space 
if every eigenvalue of $\ad_{aH} |_{\LN}$ is a positive integer for some $a >0$. 
\end{Prop}

The Ricci curvatures of solvable extensions can be calculated as follows. 

\begin{Lem} 
\label{ric-s} 
Denote the Ricci curvatures of 
$(\LS = \R H + \LN, \inner{}{}^c)$ and $(\LN, B_{\sigma})$ 
by $\ric^c$ and $\ric^{\LN}$, respectively. 
Let $\ad$ and $\ad^{\LG}$ be the adjoint operators of $\LS$ and $\LG$, 
respectively. 
We normalize $H$ to be the mean curvature vector of 
$(\LS, \inner{}{}^1 = B_{\sigma})$. 
Then, 
\begin{enumerate} 
\item[\no{1}] \ 
$\ric^c(H) = - (\tr(\ad_H)^2) / (c \cdot \tr(\ad_H^{\LG})^2) \cdot H$, 
\item[\no{2}] \ 
if $[H, \LG_0] =0$, then $\ric^c(H) = -(1/2c) H$, 
\item[\no{3}] \ 
$\ric^c(U) = \ric^{\LN}(U) - (1/c) [H,U]$ \ for every $U \in \LN$, 
\item[\no{4}] \ 
$\ric^c(U_1) = -(1/4) U_1 + [(1/2) \sum Z_k - (1/c)H, U_1]$ \ 
for every $U_1 \in \LG_1$, 
\item[\no{5}] \ 
$\ric^c(U_2) = -(1/4) U_2 + [(1/2) \sum Z_k - (1/c)H, U_2]$ \ 
for every $U_2 \in \LG_2$. 
\end{enumerate} 
\end{Lem} 

\begin{proof}
(1): \ 
Proposition \ref{wolter} leads that 
$\ric^c(H) = - (\tr(\ad_{H})^2 / |H|^2) \cdot H$. 
By the definition of the inner product, we have 
\[
|H|^2 = c \cdot B_{\sigma}(H,H) = c \cdot \tr(\ad_H^{\LG})^2 . 
\]

(2): \ 
If $[H, \LG_0] = 0$, then we have 
\begin{eqnarray*}
\tr(\ad_H^{\LG})^2 
= \tr(\ad_H^{\LG} |_{\LN})^2 + \tr(\ad_H^{\LG} |_{\LG_0})^2 + 
\tr(\ad_H^{\LG} |_{\sigma(\LN)})^2 
= 2 \tr(\ad_H)^2. 
\end{eqnarray*}
We finish the proof of (2) by substituting this for (1). 

(3): \ 
Note that the mean curvature vector of $(\LS, \inner{}{}^c)$ is $(1/c)H$. 
Thus Proposition \ref{wolter} leads that 
\[ 
\ric^c(U) = \ric^{\LN}(U) - (\tr(\ad_{(1/c)H}) / |(1/c)H|^2) \cdot [(1/c)H,U]
\] 
for $U \in \LN$. 
Recall that $\inner{H}{A}^1 = \tr(\ad_A)$ for every $A \in \R H$, 
since $H$ is the mean curvature vector of $(\LS, \inner{}{}^1)$. 
Thus one has 
\[ 
|(1/c)H|^2 = \inner{(1/c)H}{(1/c)H}^c = (1/c) \inner{H}{H}^1 
= (1/c) \tr(\ad_H) = \tr(\ad_{(1/c)H}) . 
\] 
This completes the proof of (3). 

(4) and (5): \ 
They are direct consequences of (3) and Corollary \ref{mcv-4}. 
\qed
\end{proof}

\bigskip 

Next we study the Einstein condition of the solvable extension. 
The following theorem implies that 
there is the only one natural choice of the solvable extension. 

\begin{Thm} 
\label{ric-s-1} 
Let $(\LN, \inner{}{})$ be the nilmanifold attached to 
a graded Lie algebra of type $\alpha_0$ and 
consider the above solvable extension $(\LS = \R H + \LN, \inner{}{}^c)$. 
Then we have the following$:$ 
\begin{enumerate} 
\item[\no{1}] \ 
If $(\LS, \inner{}{}^c)$ is Einstein, then $\ric^c(H) = -(1/2c) H$. 
\item[\no{2}] \ 
Assume that $\LN$ is abelian. 
Then $(\LS, \inner{}{}^c)$ is Einstein if and only if $H \parallel H_0$. 
\item[\no{3}] \ 
Assume that $\LN$ is not abelian. 
If $(\LS, \inner{}{}^c)$ is Einstein, then $c=2$ and $H \parallel H_0$. 
\end{enumerate} 
\end{Thm} 

\begin{proof}
First of all, we normalize $H$ to be the mean curvature vector of 
$(\LS, \inner{}{}^1 = B_{\sigma})$. 

(1): \ 
Assume that $(\LS, \inner{}{}^c)$ is Einstein with Einstein constant $s$. 
We have only to show $[H, \LG_0] =0$, because of Lemma \ref{ric-s} (2). 
Let $U_1 \in \LG_1$. 
Einstein condition and Lemma \ref{ric-s} (4) lead that 
\[ 
[-(1/4) Z + (1/2) \sum Z_k - (1/c)H, U_1] 
= \ric^c(U_1) = sU_1 = [sZ, U_1] . 
\] 
Thus, the element 
$U := -(1/4) Z + (1/2) \sum Z_k - (1/c)H - sZ$ 
acts trivially on $\LG_1$. 
By assumption, the gradation is of type $\alpha_0$, that is, 
$\LG$ is generated by $\LG_1$ and the grade-reversing Cartan involution. 
Therefore, 
$\sigma(U)=-U$ and $[U, \LG_1] = 0$ lead that $[U, \LG] = 0$. 
We conclude that 
\begin{eqnarray}
\label{eqn-1}
-(1/4) Z + (1/2) \sum Z_k - (1/c)H = sZ, 
\end{eqnarray}
because $\LG$ is semisimple. 
Theorem \ref{mcv-2} leads that $[Z, \LG_0]=0$ and $[Z_k, \LG_0] = 0$, 
which conclude $[H, \LG_0] = 0$. 

(2): \ 
We consider the case that $\LN$ is abelian. 
Note that, in this case, we have $H_0 = \sum Z_k = Z_1 = (1/2)Z$ 
(the last equality follows from Theorem \ref{mcv-3}). 
If $(\LS, \inner{}{}^c)$ is Einstein with Einstein constant $s$, 
the equation (\ref{eqn-1}) implies that $H$ is parallel to $H_0$. 
To show the converse, we assume that $H \parallel H_0$. 
One has $\ric^c(H_0) = -(1/2c)H_0$, 
because of $[H_0,\LG_0] =0$ and Lemma \ref{ric-s} (2). 
One also has 
\begin{eqnarray*}
\ric^c(U_1) = - (1/c) [H, U_1] = -(1/2c) [Z, U_1] = -(1/2c) U_1 , 
\end{eqnarray*}
because of $\ric^{\LN} =0$ and Lemma \ref{ric-s} (3). 
Therefore $(\LS, \inner{}{}^c)$ is Einstein 
with Einstein constant $-(1/2c)$. 

(3): \ 
We consider the case that $\LN$ is not abelian 
(that is, $\LG_2 \neq 0$) and $(\LS, \inner{}{}^c)$ is Einstein. 
For $U_1 \in \LG_1$ and $U_2 \in \LG_2$, 
Einstein condition and Lemma \ref{ric-s} (4), (5) imply that 
\begin{eqnarray*}
[(1/2) \sum Z_k - (1/c)H, U_1] = sU_1, \quad 
[(1/2) \sum Z_k - (1/c)H, U_2] = sU_2, 
\end{eqnarray*}
where the Einstein constant is $-(1/4) + s$. 
On the other hand, 
since our gradation is of type $\alpha_0$, one has $\LG_2 = [\LG_1, \LG_1]$. 
This implies that $(1/2) \sum Z_k - (1/c)H$ acts on $\LG_2$ as $2s \cdot \id$. 
Therefore $s=0$ and $H \parallel \sum Z_k = H_0$. 
Recall that we normalized $H$ to be mean curvature vector, 
which leads that $H = H_0$ and $c=2$. 
\qed 
\end{proof}

\bigskip 

In the case that $\LN$ is abelian, that is, the gradation is of first kind, 
the Einstein space $(\LS, \inner{}{}^c)$ 
is obviously isometric to the real hyperbolic space, which has been well-known. 
In the case that $\LN$ is not abelian, 
our class of solvmanifolds contains many new examples of Einstein solvmanifolds. 
We will see examples in the latter sections. 

\section{Solvmanifolds attached to graded Lie algebras}
\label{sec-naturalextension}

Theorem \ref{ric-s-1} states that 
the following solvable extension is natural: 

\begin{Def}
{\rm 
Let $(\LN, B_{\sigma})$ be the nilmanifold attached to a graded Lie algebra 
and $H_0$ the mean curvature vector of $(\LA + \LN, B_{\sigma})$. 
We call $(\LS = \R H_0 + \LN, \inner{}{}^2)$ the 
{\it solvmanifold attached to a graded Lie algebra}. 
}
\end{Def}

Note that the solvmanifold attached to GLA is a Carnot space, 
because of Theorem \ref{mcv-2}. 
On the Ricci curvatures of these solvmanifolds, 
Lemma \ref{ric-s} immediately leads the following: 

\begin{Thm}
\label{ricciofoursolv}
The Ricci curvature $\ric$ of a solvmanifold 
attached to a graded Lie algebra, $(\LS = \R H_0 + \LN, \inner{}{}^2)$, 
satisfies that 
\begin{enumerate}
\item[\rm{(1)}] 
$\ric = -(1/4) \id$ on $\R H_0 + \LG_1 + \LG_2$, 
\item[\rm{(2)}] 
$\ric(U) = \ric^{\LN}(U) - (1/2) [H_0,U]$ \ for every $U \in \LN$. 
\end{enumerate}
\end{Thm}

We here mention some relations between our solvmanifolds 
and symmetric spaces of noncompact type. 

Let us consider $\LG = \sum_{k=-\nu}^{\nu} \LG_k$ the longest gradation. 
Here, we call the gradation {\it longest} if 
the characteristic element $Z$ is a regular element. 
In this case, 
$\LN = \sum_{k>0} \LG_k = \sum_{\alpha \in \Delta^+} \LG_{\alpha}$ 
coincides with the nilpotent part of the Iwasawa decomposition of $\LG$. 
Furthermore, the solvmanifold $(\LA + \LN, \inner{}{})$, where 
$\inner{}{} := 2 B_{\sigma} |_{\LA \times \LA} + B_{\sigma} |_{\LN \times \LN}$, 
is nothing but the symmetric space of noncompact type, $G/K$, associated from $\LG$. 
This is a standard Einstein solvmanifold, 
whose rank one reduction coincides with our solvable extension 
$(\R H_0 + \LN, \inner{}{}^2)$. 
Thus, the natural solvable extension is Einstein if the gradation is longest. 

Let us consider the gradation $\LG = \sum_{k=-\nu}^{\nu} \LG_k$, 
which is not necessarily longest. 
In this case we have 
$\LN = \sum_{k>0} \LG_k \subset \sum_{\alpha \in \Delta^+} \LG_{\alpha}$. 
Therefore, our solvable Lie algebra $\R H_0 + \LN$ is a subalgebra 
of the solvable part of the Iwasawa decomposition. 
Thus one can see that our solvmanifold $(\R H_0 + \LN, \inner{}{}^2)$ is 
a Riemannian submanifold in a symmetric space of noncompact type. 

Because of the above observation, 
we have to investigate when our solvmanifold is totally geodesic, 
or the rank one reduction of a totally geodesic submanifold. 
As a result, this can only happen if the gradation is longest. 

\begin{Prop}
Let $\LG = \sum \LG_k$ be a gradation and $\LN := \sum_{k>0} \LG_k$. 
If $\LN$ is the nilradical of a totally geodesic submanifold, 
then the gradation is longest. 
\end{Prop}

\proof
Assume that $\LA' + \LN$ gives a totally geodesic submanifold 
and $\R H_0 + \LN$ is the rank one reduction of it. 
Denote by $\LG = \LK + \LP$ the Cartan decomposition and 
$\pi : \LG \rightarrow \LP$ the natural projection. 
One knows that $V := \pi(\LA' + \LN)$ is a Lie triple system. 
If we denote by $Z$ the characteristic element of our gradation, 
we have 
\[
V = \LA' + \sum_{\alpha(Z)>0} \LP_{\alpha} . 
\]
The totally geodesic submanifold coincides with the orbit of a Lie group 
with Lie algebra $[V,V] + V$. 

We here claim that $[V,V] + V = \LG$. 
Let us take $\alpha$ such that $\alpha(Z)>0$. 
Theorem \ref{mcv-2} leads $\alpha(H_0)>0$, and thus one has 
$[V,V] \supset [H_0, \LP_{\alpha}] = \LK_{\alpha}$. 
This implies $\LG_{-1} + \LG_1 \subset [V,V] + V$. 
Since we assumed that the gradation is of type $\alpha_0$, 
we conclude the claim. 

Now it is easy to see from $\pi(\LA' + \LN) = \LP$ that 
the totally geodesic submanifold is the entire space, 
and $Z$ is regular. 
\qed

\section{Two-step case} 
\label{sec-two}

In this section we see that 
every solvmanifold attached to a graded Lie algebra 
of second kind is Einstein. 
Most of those Einstein solvmanifolds have 
the eigenvalue type $(1<2 \ ; \ d_1,d_2)$, 
but some of them have another eigenvalue type. 
We refer \cite{K} for the explicit classification list 
of second kind graded Lie algebras. 

\begin{Thm} 
\label{secondkind}
Every solvmanifold $(\LS, \inner{}{})$ attached to 
a graded Lie algebra of second kind is Einstein. 
\end{Thm} 

\proof 
Since $\LS = \R H_0 + \LG_1 + \LG_2$, 
this is a direct consequence of Theorem \ref{ricciofoursolv}. 
\qed 

\bigskip 

We study their eigenvalue types. 
For a graded Lie algebra of second kind, $\LG = \sum_{k=-2}^2 \LG_k$, 
the characteristic element satisfies either $Z = H^i$ or $Z = H^i + H^j$. 
Recall that $\{ H^i \}$ is the dual basis of the simple roots. 
If $Z = H^i$, Theorem \ref{mcv-2} leads that $Z \parallel H_0$. 
In this case the eigenvalue type is $(1<2 \ ; \ d_1,d_2)$, 
by the definition of the characteristic element. 
If $Z = H^i + H^j$, Theorem \ref{mcv-2} states only that 
$H_0 = a H^i + b H^j$. 
In this case unusual eigenvalue type may occur. 

Let us consider the root system is of $A_r$-type. 
In this case $Z = H^i + H^j$ gives a second kind gradation. 
One can calculate the mean curvature vector directly, 
\[ 
H_0 = Z_1 + Z_2 = j H^i + (r+1-i) H^j , 
\] 
up to scalar. 
If $i+j = r+1$ then the eigenvalue type is $(1<2 \ ; \ d_1, d_2)$. 
If $i+j < r+1$, the eigenvalue type is 
$(j < r+1-i < r+1-i+j \ ; \ d_1, d_2, d_3)$. 
These arguments imply that 

\begin{Cor} 
For every $0< \mu_1 < \mu_2$, 
there exists a standard Einstein solvmanifold $(\LS, \inner{}{})$
with two-step nilpotent $\LN$,
whose eigenvalue type is
$(\mu_1 < \mu_2 < \mu_1 + \mu_2 \ ; \ d_1, d_2, d_3)$.
\end{Cor} 

We see a simple example here. 
Consider the nilpotent subalgebra of $\sll(n+m+l,\R)$, 
\[ 
\LN := \left\{ \left( \begin{array}{ccc} 
0 & X & Z \\ 0 & 0 & Y \\ 0 & 0 & 0 
\end{array} \right) \mid 
X \in M_{n,m}(\R), Y \in M_{m,l}(\R), Z \in M_{n,l}(\R) 
\right\} . 
\] 
Let us consider the gradation of $\sll(n+m+l,\R)$ of second kind, 
given by the above block decomposition. 
Then $(\LN, B_{\sigma})$ coincides with the nilmanifold attached to 
this gradation, 
therefore it has an Einstein solvable extension. 

In this case, the root system is of $A_{n+m+l-1}$-type, 
and direct calculations lead that 
\begin{eqnarray*}
H^n & = & \frac{1}{n+m+l} \left( 
\begin{array}{ccc}
(m+l) \cdot 1_n & & \\ 
& -n \cdot 1_m & \\ 
& & -n \cdot 1_l 
\end{array}
\right) , \\ 
H^m & = & \frac{1}{n+m+l} \left( 
\begin{array}{ccc}
l \cdot 1_n & & \\ 
& l \cdot 1_m & \\ 
& & -(n+m) \cdot 1_l 
\end{array}
\right) . 
\end{eqnarray*}
The characteristic element is $Z=H^n+H^m$. 
The solvmanifold $\R Z + \LN$ is a three-step Carnot space, 
but this is not what we want. 
Our Einstein solvable extension is $\R H_0 + \LN$, 
where $H_0$ is parallel to $(n+m)H^n + (m+l)H^{n+m}$. 
This is $(n+2m+l+1)$-step Carnot space 
if $n+m$ and $m+l$ have no common devisor, 
whereas $\LN$ is two-step nilpotent. 

We also note that, in case $n=l=1$, this nilmanifold is nothing but the 
Heisenberg Lie algebra, 
and the constructed Einstein solvmanifold is isometric to 
the complex hyperbolic space. 
In case $n=l$, the eigenvalue type of the Einstein solvmanifold 
is $(1<2 ; nm+ml, nl)$. 
In case $n \neq l$, $\ad_{H_0}$ has three distinct eigenvalues. 

\section{Examples of Einstein spaces with higher steps}
\label{higher-step}

In this section we investigate the solvmanifolds 
attached to graded Lie algebras of third kind and of fourth kind. 

\begin{Prop} 
\label{3-step} 
Assume a graded Lie algebra is of third kind and 
each $Z_k$ is parallel to $Z$. 
Then the attached solvmanifold $(\LS, \inner{}{})$ is Einstein 
if and only if $\dim \LG_1 = 2\dim \LG_2$. 
\end{Prop} 

\proof 
Recall that Theorem \ref{ricciofoursolv} states 
$\ric = -(1/4) \cdot \id$ on $\R H_0 + \LG_1 + \LG_2$. 
Therefore $(\LS, \inner{}{})$ is Einstein if and only if 
$\ric(U_3) = -(1/4) U_3$ for every $U_3 \in \LG_3$. 
One can calculate this as follows: 
\begin{eqnarray*} 
\ric^{\LN}(U_3) & = & (1/4)[Z_1,U_3] + (1/4)[Z_2, U_3] , \\
\ric(U_3) & = & \ric^{\LN}(U_3) - (1/2) [Z_1 + Z_2 + Z_3,U_3] \\
& = & - (1/4) [Z_1,U_3] - (1/4)[Z_2,U_3] -(1/2)[Z_3,U_3] \\
& = & -(1/4) ((\dim \LG_1 + 2 \dim \LG_2 + 6 \dim \LG_3)/B_{\sigma}(Z,Z))
\cdot [Z, U_3] \\ 
& = & -(1/4) (3\dim \LG_1 + 6 \dim \LG_2 + 18 \dim \LG_3)/B_{\sigma}(Z,Z))
\cdot U_3. 
\end{eqnarray*} 
Thus, $\ric(U_3) = -(1/4) U_3$ if and only if 
\[ 3 \dim \LG_1 + 6 \dim \LG_2 + 18 \dim \LG_3 = B_{\sigma}(Z,Z)
= 2(\dim \LG_1 + 4\dim \LG_2 + 9 \dim \LG_3), \] 
that is, $\dim \LG_1 = 2\dim \LG_2$. 
\qed 

\bigskip 

Next we state the similar result for four-step case. 
Proof is also similar. 

\begin{Prop} 
\label{4-step} 
Assume a graded Lie algebra is of fourth kind and 
each $Z_k$ is parallel to $Z$. 
Then the attached solvmanifold $(\LS, \inner{}{})$ is Einstein 
if and only if $\dim \LG_1 = 3 \dim \LG_3$ and 
$\dim \LG_1 + 4 \dim \LG_4 = 2 \dim \LG_2$. 
\end{Prop} 

These propositions provide examples of Einstein solvmanifolds 
whose nilradicals are of three-step or four-step. 

\begin{Thm} 
\label{listof3and4}
The solvmanifolds attached to the following graded Lie algebras 
are Einstein. 
The eigenvalue types of the constructed Einstein solvmanifolds 
are as follows. 
\end{Thm} 

\[ \begin{array}{ccccl}
\quad \LG \quad & \quad \Delta \quad & \quad Z \quad && \mbox{eigenvalue
type} \\ \hline \hline
\LG_{2(2)} & G_2 & H^1 && (1<2<3 \ ; \ 2,1,2) \\
\LG_2^{\C} & G_2 & H^1 && (1<2<3 \ ; \ 4,2,4) \\ \hline
\LF_{4(4)} & F_4 & H^2 && (1<2<3 \ ; \ 12,6,2) \\
\LF_4^{\C} & F_4 & H^2 && (1<2<3 \ ; \ 24,12,4) \\
\LE_{6(2)} & F_4 & H^2 && (1<2<3 \ ; \ 18,9,2) \\
\LE_{7(-5)} & F_4 & H^2 && (1<2<3 \ ; \ 30,15,2) \\
\LE_{8(-24)} & F_4 & H^2 && (1<2<3 \ ; \ 54,27,2) \\ \hline
\LE_{6(6)} & E_6 & H^4 && (1<2<3 \ ; \ 18,9,2) \\
\LE_6^{\C} & E_6 & H^4 && (1<2<3 \ ; \ 36,18,4) \\ \hline
\LF_{4(4)} & F_4 & H^3 && (1<2<3<4 \ ; \ 6,9,2,3) \\
\LF_4^{\C} & F_4 & H^3 && (1<2<3<4 \ ; \ 12,18,4,6) \\
\LE_{6(2)} & F_4 & H^3 && (1<2<3<4 \ ; \ 12,12,4,3) \\
\LE_{7(-5)} & F_4 & H^3 && (1<2<3<4 \ ; \ 24,18,8,3) \\
\LE_{8(-24)} & F_4 & H^3 && (1<2<3<4 \ ; \ 48,30,16,3) \\ \hline
\end{array} \] 

\begin{proof}
One has to check the relations of the dimensions. 
One can calculate $\dim \LG_k$  by counting roots and their multiplicities 
(one can see the table of multiplicities in \cite{T2}). 
\qed
\end{proof}

\section{The case of low rank}
\label{sec-lowrank}

In this section we mention the following problem: \ 
for a given semisimple Lie algebra $\LG$, 
is every attached solvmanifold Einstein? \ 
The answer is yes, 
if the rank of the associated symmetric space $(\LG, \LK)$ is low. 

If the rank of $(\LG, \LK)$ is one, it is trivial. 
In this case, $\LG$ has only one gradation, which is longest, 
and the attached solvmanifold is the entire symmetric space. 

\begin{Prop}
Assume that the rank of $(\LG, \LK)$ is two, that is, 
$\LG$ is one of the following$:$ 
\[
\mathrm{sl}_3(\R), \ 
\mathrm{sl}_3(\C), \ 
\mathrm{su}^{\ast}(6), \ 
\mathrm{e}_{6(-26)}, \ 
\mathrm{so}(5, \C), \ 
\mathrm{so}(2, n+2), \ 
\LG_{2(2)}, \ 
\LG_2^{\C}. 
\]
Then, for every gradation of $\LG$, the attached solvmanifold is Einstein. 
\end{Prop}

\proof 
Let $\{ \alpha_1, \alpha_2 \}$ be the set of simple roots. 
There are three gradations, given by $H^1$, $H^2$ and $H^1+H^2$. 
If the characteristic element is $H^1 + H^2$, 
the gradation is longest and thus the attached solvmanifold is Einstein. 
If the characteristic element is $H^1$ or $H^2$, 
the gradation is of first kind, second kind or third kind. 
The solvmanifold attached to first kind or second kind gradation 
is Einstein (Theorems \ref{ricciofoursolv}, \ref{secondkind}). 
The third kind gradation is obtained by $G_2$-type root system, 
whose attached solvmanifolds are Einstein by Theorem \ref{listof3and4}. 
\qed

\bigskip 

Similarly, we have the following: 

\begin{Prop}
If the root system of $(\LG, \LK)$ is of $A_3$-type, 
that is, 
$\LG = \mathrm{sl}_4(\R)$, $\mathrm{sl}_4(\C)$ or $\mathrm{su}^{\ast}(8)$, 
every attached solvmanifold is Einstein. 
\end{Prop}

\section{Rank one reduction of noncompact symmetric spaces}
\label{sec-symm} 

In this section we study 
the mean curvature vectors of noncompact symmetric spaces. 
For the graded Lie algebra with characteristic element 
$Z := H^1 + \cdots + H^r$, 
one knows that $\LN$ is the nilpotent part of the Iwasawa decomposition 
and the solvmanifold $(\LA + \LN, \inner{}{})$ is a symmetric space. 

\begin{Prop}
The mean curvature vector of the above solvmanifold 
$(\LA + \LN, \inner{}{})$ is $H_0 = h_1H^1 + \cdots + h_rH^r$, 
where $h_i = (\dim \LG_{\alpha_i} + 4 \dim \LG_{2\alpha_i}) 
B_{\sigma} (H_{\alpha_i},H_{\alpha_i})$. 
\end{Prop}

\proof 
We apply Theorem \ref{mcv-2} (2). 
In this case $\Delta^{i} := \{ \beta \mid \beta(Z) = \beta(H^i) >0 \}$ 
consists of $\alpha_i$ and $2\alpha_i$ (if exists). 
Therefore, 
\[ h_i = (\dim \LG_{\alpha_i}) B_{\sigma} (H_{\alpha_i},H_{\alpha_i}) 
+ (\dim \LG_{2\alpha_i}) B_{\sigma} (H_{2\alpha_i},H_{2\alpha_i}) . \]
Proposition is completed by $H_{2\alpha_i} = 2 H_{\alpha_i}$. 
\qedpar 

This Proposition implies that, 
if every root has the same length, then $H_0$ is parallel to $Z$. 
It is remarkable that 
$H_0$ may not parallel to $Z$ even for symmetric spaces. 
Furthermore, as we saw in Propositions \ref{3-step} and \ref{4-step}, 
the Einstein condition requires some equalities of dimensions. 
It seem to be mysterious that symmetric spaces automatically 
satisfy these conditions.

\end{document}